\documentclass[11pt,leqno]{article}
\usepackage{amsmath}
\usepackage{amsthm}

\numberwithin{equation}{section}
\let\sect=\section

\newtheorem{theorem}{Theorem}[section]
\newtheorem{corollary}[theorem]{Corollary}
\newtheorem{lemma}[theorem]{Lemma}
\newtheorem{proposition}[theorem]{Proposition}
\newtheorem{claim}[theorem]{Claim}
\newtheorem{example}[theorem]{\sl Example}

\theoremstyle{definition}
\newtheorem{remark}[theorem]{Remark}

\newcommand{\EE}{{\bf  E}}

\newcommand{\RR}{{\bf  R}}

\newcommand{\CC}{{\bf  C}}

\newcommand{\PP}{{\bf  P}}

\newcommand{\Ct}{{\tilde{C}}}

\newcommand{\Wt}{{\widetilde{W}}}

\newcommand{\Zt}{{\tilde{Z}}}

\newcommand{\Cc}{{\cal C}}

\newcommand{\Fc}{{\cal F}}

\newcommand{\Mc}{{\cal M}}

\newcommand{\begp}{\begin{proposition}}
\newcommand{\enp}{\end{proposition}}
\newcommand{\begt}{\begin{theorem}}
\newcommand{\ent}{\end{theorem}}
\newcommand{\begl}{\begin{lemma}}
\newcommand{\enl}{\end{lemma}}
\newcommand{\begc}{\begin{corollary}}
\newcommand{\enc}{\end{corollary}}
\newcommand{\begcl}{\begin{claim}}
\newcommand{\encl}{\end{claim}}
\newcommand{\begr}{\begin{remark}}
\newcommand{\enr}{\end{remark}}
\newcommand{\begal}{\begin{algorithm}}
\newcommand{\enal}{\end{algorithm}}
\newcommand{\begd}{\begin{definition}}
\newcommand{\enf}{\end{definition}}
\newcommand{\begx}{\begin{example}}
\newcommand{\enx}{\end{example}}
\newcommand{\bega}{\begin{array}}
\newcommand{\ena}{\end{array}}

\def\rompar(#1){\textup(#1\textup)}    

\newcommand{\refT}[1]{Theorem~\ref{#1}}
\newcommand{\refC}[1]{Corollary~\ref{#1}}
\newcommand{\refL}[1]{Lemma~\ref{#1}}

\newcommand{\refP}[1]{Proposition~\ref{#1}}

\newcommand\nopf{\qed}   

\newcommand\Roesler{R\"{o}sler}

\begin{document}

\setcounter{page}{0}
\thispagestyle{empty}

\begin{center}
{\Large \bf A Characterization of the Set of Fixed Points of the Quicksort Transformation
\\ } 
\normalsize
 
\vspace{8ex}
{\sc James Allen Fill\footnotemark} \\
\vspace{.1in}
Department of Mathematical Sciences \\ 
\vspace{.1in}
The Johns Hopkins University \\
\vspace{.1in}
{\tt jimfill@jhu.edu} and {\tt http://www.mts.jhu.edu/\~{}fill/} \\
\vspace{.2in}
{\sc and} \\
{\sc Svante Janson}\\ 
\vspace{.1in}
Department of Mathematics \\ 
\vspace{.1in}
Uppsala University \\
\vspace{.1in}
{\tt svante.janson@math.uu.se} and {\tt http://www.math.uu.se/\~{}svante/} \\
\end{center}
\vspace{3ex}
 
\begin{center}
{\sl ABSTRACT} \\
\end{center}
\vspace{.05in}

\begin{small}
The limiting distribution~$\mu$ of the normalized number of key comparisons required by
the {\tt Quicksort} sorting algorithm is known to be the unique fixed point
of a certain distributional transformation~$T$---unique, that is, subject to the constraints
of zero mean and finite variance.  We show that a distribution is a fixed point of~$T$
if and only if it is the convolution of~$\mu$ with a Cauchy distribution of arbitrary center
and scale.  In particular, therefore, $\mu$ is the unique fixed point of~$T$ having zero mean.
\smallskip
\par\noindent
{\em AMS\/} 2000 {\em subject classifications.\/}  Primary 68W40;
secondary 60E05, 60E10, 68P10.
\medskip
\par\noindent
{\em Key words and phrases.\/} {\tt Quicksort}, fixed point, characteristic function, smoothing
transformation, domain of attraction, coupling, integral equation.
\medskip
\par\noindent
\emph{Date.} May~22, 2000.
\end{small}

\footnotetext[1]{Research supported by NSF grant DMS--9803780,
and by the Acheson J.~Duncan Fund for the Advancement of Research in Statistics.}

\newpage
\addtolength{\topmargin}{+0.5in}

\sect{Introduction, motivation, and summary}
\label{intro}

Let~$\Mc$ denote the class of all probability measures (distributions) on the real line.
This paper concerns the transformation~$T$ defined on~$\Mc$ by letting~$T \nu$
be the distribution of
$$
U Z + (1 - U) Z^* + g(U),
$$
where~$U$, $Z$, and~$Z^*$ are independent, with $Z \sim \nu$, $Z^* \sim \nu$,
and $U \sim \mbox{unif}(0, 1)$, and where
\begin{equation}
\label{gdef}
g(u) := 2 u \ln u + 2 (1 - u) \ln (1 - u) + 1.
\end{equation}
Of course, $T$ can be regarded as a transformation on the class of
characteristic functions~$\psi$ of elements of~$\Mc$.  With this interpretation,
$T$ takes the form
$$
(T \psi)(t) = \EE \big[ \psi(U t)\,\psi((1 - U) t)\,\exp[i t g(U)] \big]
            = \int^1_{u = 0}\!\psi(u t)\,\psi((1 - u) t)\,e^{i t g(u)}\,du,\ \ t \in \RR.
$$

It is well known~\cite{Roesler} that (i)~among distributions with zero mean and finite
variance, $T$ has a unique fixed point, call it~$\mu$; and
(ii)~if $C_n$ denotes the random number of key comparisons required by the algorithm
{\tt Quicksort} to sort a file of~$n$ records, then the distribution of
$(C_n - \EE C_n) / n$ converges weakly to~$\mu$.

There are other fixed points.  For example, it has been noted frequently in the
literature that the location family generated by~$\mu$ is a family of fixed points.
But there are many more fixed points, as we now describe.
Define the Cauchy($m, \sigma$) distribution (where $m \in \RR$ and $\sigma \geq 0$)
to be the distribution of $m + \sigma C$, where~$C$ has the standard Cauchy distribution
with density $x \mapsto [\pi (1 + x^2)]^{-1}$, $x \in \RR$;
equivalently, Cauchy($m, \sigma$) is the distribution with characteristic
function $e^{i m t - \sigma |t|}$.  [In particular,
the Cauchy($m, 0$) distribution is unit mass at~$m$.]
Now let~$\Fc$ denote the class of all fixed points of~$T$, and
let~$\Cc$ denote the class of convolutions of~$\mu$ with a Cauchy distribution.
Using characteristic functions it is easy to check
that $\Cc \subseteq \Fc$, and that all of the distributions in~$\Cc$ are distinct.
In this paper we will prove that, conversely, $\Fc \subseteq \Cc$, and thereby
establish the following main result.

\begin{theorem}
\label{T:char}
The class $\Fc$ equals~$\Cc$.  That is, a measure~$\nu$ is a fixed point of the
{\tt Quicksort} transformation~$T$ if and only if it is the convolution of the
limiting {\tt Quicksort} measure~$\mu$ with a Cauchy distribution of arbitrary
center~$m$ and scale~$\sigma$.  In particular, $\Fc$ is in one-to-one correspondence
with the set $\{ (m, \sigma): m \in \RR,\ \sigma \geq 0 \}$.
\end{theorem}

The following corollary is immediate and strengthens \Roesler's~\cite{Roesler}
characterization of~$\mu$ as the unique element of~$\Fc$ having zero mean
and finite variance.

\begin{corollary}
\label{C:zeromean}
The limiting {\tt Quicksort} measure~$\mu$ is the unique fixed point of the
{\tt Quicksort} transformation~$T$ having finite expectation equal to~$0$.
\nopf
\end{corollary}

The present paper can be motivated in two ways.  First, the authors are writing
a series of papers refining and extending \Roesler's~\cite{Roesler} probabilistic
analysis of {\tt Quicksort}.  No closed-form expressions are known for any of the
standard functionals (e.g.,\ characteristic function, distribution function,
density function) associated with~$\mu$; information to be obtained about~$\mu$
must be read from the fixed-point identity it satisfies.  We were curious as to
what extent additional known information about~$\mu$, such as the fact that it has
everywhere finite moment generating function, must be brought to bear.  As one example,
it is believed that the continuous Lebesgue density~$f$
(treated in~\cite{quick_density}) for~$\mu$ decays at least exponentially quickly
to~$0$ at $\pm \infty$, cf.~\cite{KnSz}.  But we now know from \refT{T:char}
that there there can be no proof for this conjecture
that solely makes use of the information that $\mu \in \Fc$.

Second, we view the present paper as a pilot study of fixed points for a {\em class\/} of 
distributional transformations on the line.  In the more general setting, we would
be given (the joint distribution of) a sequence $(A_i: i \geq 0)$ of
random variables and would define a transformation~$T$ on~$\Mc$ by letting~$T \nu$
be the distribution of $A_0 + \sum_{i = 1}^{\infty} A_i Z_i$, where $Z_1, Z_2, \ldots$
are independent random variables with distribution~$\nu$.  [To ensure well-definedness,
one might (for example) require that (almost surely) $A_i \neq 0$ for only finitely
many values of~$i$.]  For probability measures~$\nu$ on $[0, \infty)$, rather than
on~$\RR$, and with the additional restrictions that $A_0 = 0$ and $A_i \geq 0$
for all $i \geq 1$, such transformations are called {\em generalized smoothing
transformations\/}.  These have been thoroughly studied by Durrett and Liggett~\cite{DL},
Guivarc'h~\cite{Guiv}, and Liu~\cite{Liu}, and by other authors; consult the three
papers we have cited here for further bibliographic references.
Generalized smoothing transformations have applications to interacting particle
systems, branching processes and branching random walk, random set constructions, and
statistical turbulence theory.
The arguments used to characterize the set of fixed points for generalized smoothing
transformations make heavy use of Laplace transforms; unfortunately, these arguments do not
carry over readily to distributions on the line.  Other authors
(see, e.g.,~\cite{Roesler92}~\cite{Roesler99}~\cite{RR}) have treated fixed
points of transformations of measures~$\nu$ on the whole line as discussed above, but not
without finiteness conditions on the moments of~$\nu$.

We now outline our proof of \refT{T:char}.  Let~$\psi$ be the characteristic function 
of a given $\nu \in \Fc$,
and let $r(t) := \psi(t) - 1$, $t \in \RR$.
In Section~\ref{inteq} we establish and solve (in a certain sense)
an integral equation satisfied by~$r$.
In Section~\ref{asy} we then use the method of successive substitutions to
derive asymptotic information about~$r(t)$ as $t \downarrow 0$, showing first that $r(t) =
O(t^{2 / 3})$, next that $r(t) = \beta t + o(t)$ for some
$\beta = -\sigma + i m \in \CC$ with $\sigma \geq 0$,
and finally that $r(t) = \beta t + O(t^2)$ .
In Section~\ref{main} we use this information to argue
that there exist random variables $Z_1 \sim \nu$,
$Z_2 \sim \mu$, and $C \sim \mbox{Cauchy}(m, \sigma)$
such that $Z_1 = Z_2 + C$.  We finish the proof by showing
that one can take~$Z_2$ and~$C$ to be independent, whence $\nu \in \Cc$.

\sect{An integral equation}
\label{inteq}

Let~$\psi$ denote the characteristic function of a given $\nu \in \Fc$.
Since $\psi(-t) \equiv \overline{\psi(t)}$, we shall only need to consider $\psi(t)$
for $t \geq 0$.
For notational convenience, define
$$
r(t) := \psi(t) - 1,\ \ t \geq 0.
$$
Rearranging the fixed-point integral equation $(T \psi)(t) \equiv \psi(t)$,
we obtain the following result.

\begin{lemma}
\label{L:rinteq}
The function~$r$ satisfies the integral equation
$$
r(t) = 2 \int_{u = 0}^1\!r(u t)\,du + b(t),\ \ t \geq 0,
$$
where
\begin{equation}
\label{bdef}
b(t) := \int_{u = 0}^1\!r(u t)\,r((1 - u) t)\,du
          + i t \int_{u = 0}^1\![\psi(u t)\,\psi((1 - u) t) - 1]\,g(u)\,du + a(t)
\end{equation}
with
\begin{equation}
\label{adef}
\left| a(t)
  :=  \int_{u = 0}^1\!\psi(u t)\,\psi((1 - u) t)\,[e^{i t g(u)} - 1 - i t g(u)]\,du \right|
     \leq  \mbox{$\frac{1}{2} \EE g^2(U) t^2 = (\frac{7}{6} - \frac{1}{9} \pi^2) t^2$.}
\end{equation}
\nopf
\end{lemma}

Note that~$r$ and~$b$ are continuous on $[0, \infty)$, with $r(0) = 0 = b(0)$.
Regarding~$b$ as ``known'', the integral equation in \refL{L:rinteq} is easily ``solved''
for~$r$:

\begin{proposition}
\label{P:rsol}
For some constant $c \in \CC$, we have
$$
\frac{r(t)}{t} = c - 2 \int_{v = t}^1\!\frac{b(v)}{v^2}\,dv + \frac{b(t)}{t},\ \ t > 0.
$$
\end{proposition}

\begin{proof}
Setting $h(t) := t [r(t) - b(t)]$, \refL{L:rinteq} implies
$$
h(t) = 2 \int_{v = 0}^t\!\left[ \frac{h(v)}{v} +  b(v) \right]\,dv,\ \ t > 0.
$$
Thus~$h$ is continuously differentiable on $(0, \infty)$ and satisfies the differential
equation
$$
h'(t) = \mbox{$\frac{2}{t} h(t) + 2 b(t)$}
$$
there.  This is an easy differential equation to solve for~$h$, and we find that
$$
h(t) = c t^2 - 2 t^2 \int_{v = t}^1\!\frac{b(v)}{v^2}\,dv,\ \ t > 0,
$$
for some $c \in \CC$.  After rearrangement, the proposition is proved.  
\end{proof}

\sect{Behavior of~$r$ near~$0$}
\label{asy}

We now proceed in stages, using \refP{P:rsol} as our basic tool, to get ever more information
about the behavior of~$r$ (especially near~$0$).

\begin{lemma}
\label{L:2/3}
Let $\psi \equiv 1 + r$ denote the characteristic function of a given $\nu \in \Fc$.
Then there exists a constant $C < \infty$ such that
$$
\mbox{\rm $|r(t)| \leq C t^{2 / 3}$ for all $t \geq 0$.}
$$
\end{lemma}

\begin{proof}
Let
$$
M(t) := \max\{ |r(s)|: 0 \leq s \leq t \} \leq 2,\ \ t \geq 0.
$$
From~\eqref{bdef} and~\eqref{adef}, we see immediately that, for $0 < t \leq 1$,
$$
|b(t)| \leq M^2(t) + O(t).
$$
Therefore, for $0 < t < 1$, \refP{P:rsol} yields
$$
|r(t)| \leq M^2(t) + 2 t \int_{v = t}^1\!\frac{M^2(v)}{v^2}\,dv + \epsilon(t)
         = M^2(t) + 2 \int_{u = t}^1\!M^2(t / u)\,du + \epsilon(t),
$$
where
$$
\epsilon(t) =  \mbox{$O \left( t \log \left( \frac{1}{t} \right) \right) + O(t)
= O(t^{2 / 3})$.}
$$
Consequently, again for $0 < t < 1$ (but then trivially for all $t \geq 0$),
$$
M(t) \leq M^2(t) + 2 \int_{u = 0}^1\!M^2(t / u)\,du + O(t^{2 / 3}).
$$

Fix $0 < a < 1$; later in the proof we shall see that $a = 1 / 8$ suffices for our purposes.
Since $M(t) \to 0$ as $t \to 0$, we can choose $t_0 > 0$ such that $M(t_0) \leq a$.
Then, for $0 \leq t \leq t_0$,
\begin{eqnarray*}
M(t) &\leq& M^2(t) + 2 \int_{u = 0}^{t / t_0}\!M^2(t / u)\,du 
              + 2 \int_{u = t / t_0}^1\!M^2(t / u)\,du + O(t^{2 / 3}) \\
     &\leq& a M(t) + 8 \frac{t}{t_0}
              + 2 a \int_{u = t / t_0}^1\!M(t / u)\,du + O(t^{2 / 3})
\end{eqnarray*}
and thus
$$
M(t) \leq \frac{2 a}{1 - a} \int_{u = 0}^1\!M(t / u)\,du + O(t^{2 / 3}).
$$
Since~$M$ is bounded, this is trivially true also for $t > t_0$.
Summarizing, for some constant $\Ct < \infty$ we have, with $U \sim \mbox{unif}(0, 1)$, 
\begin{equation}
\label{Mbd}
M(t) \leq \frac{2 a}{1 - a} \EE M(t / U) + \Ct t^{2 / 3},\ \ t \geq 0.
\end{equation}

Now fix the value of~$a$ to be any number in $(0, 1 / 7)$, say $a = 1 / 8$.
Then a straightforward induction [substituting~\eqref{Mbdn} into~\eqref{Mbd} for the
induction step] shows that for any nonnegative integer~$n$ we have, for all $t \geq 0$,
\begin{equation}
\label{Mbdn}
M(t) \leq \left( \frac{2a}{1 - a} \right)^n \EE\,M \left( \frac{t}{U_1 \ldots U_n} \right)
            + \frac{1 - a}{1 - 7 a} \Ct t^{2 / 3}.
\end{equation}
Recalling that~$M$ is bounded and letting $n \to \infty$, we obtain the desired conclusion,
with $C := \frac{1 - a}{1 - 7 a} \Ct$.
\end{proof}

\begin{lemma}
\label{L:1}
Let $\psi \equiv 1 + r$ denote the characteristic function of a given $\nu \in \Fc$,
and define~$b$ by~\eqref{bdef}.
Then
$$
\mbox{\rm $r(t) = (c - 2 J) t + o(t)$ as $t \downarrow 0$,} 
$$
where~$J$ is the absolutely convergent integral 
\begin{equation}
\label{Jdef}
J := \int_{v = 0}^1\!\frac{b(v)}{v^2}\,dv.
\end{equation}
\end{lemma}

\begin{proof}
Combining~\eqref{bdef}--\eqref{adef} and \refL{L:2/3}, we obtain
$$
|b(t)| \leq O(t^{4 / 3}) + O(t^{1 + (2 / 3)}) + O(t^2) = O(t^{4 / 3}).
$$
Thus the integral~$J$ converges absolutely,
and from \refP{P:rsol} we obtain the desired conclusion about~$r$. 
\end{proof}

\refL{L:1} is all we will need in the next section, but the following refinement
follows readily and takes us as far as we can go with the
method of successive substitutions.

\begin{corollary}
\label{C:expand}
Let~$\psi$ denote the characteristic function of a given $\nu \in \Fc$.
Then there exists a constant $\beta = i m - \sigma \in \CC$ with $\sigma \geq 0$
such that
$$
\mbox{\rm $\psi(t) = 1 + \beta t + O(t^2)$ as $t \downarrow 0$.}
$$
\end{corollary}

\begin{proof}
Combining~\eqref{bdef}--\eqref{adef} and \refL{L:1}, we readily obtain
$b(t) = O(t^2)$.  Therefore, by \refP{P:rsol},
$$
\psi(t) - 1 = r(t) = (c - 2 J) t + 2 t \int_{v = 0}^t\!\frac{b(v)}{v^2}\,dv + b(t)
            = \beta t + O(t^2),
$$
with $\beta = i m - \sigma := c - 2 J$.  Since $|\psi(t)| \leq 1$ for all~$t$,
we must have $\sigma \geq 0$.
\end{proof}

\sect{Proof of the main theorem}
\label{main}

\subsection{Further preliminaries}
\label{prelim}

In Sections~\ref{prelim}--\ref{finish} we complete the proof of our main \refT{T:char}.
To do this, we begin with a key result that any characteristic function
with expansion as in \refC{C:expand} [more generally, we allow the remainder
term there to be simply $o(t)$] is in the domain of attraction of (iterates of)
the ``homogeneous'' analogue~$T_0$ of~$T$.  (Here $\implies$ denotes weak convergence
of probability measures.)

\begin{theorem}
\label{T:attract}
Let~$\psi$ be any characteristic function satisfying
\begin{equation}
\label{linear}
\mbox{\rm $\psi(t) = 1 + \beta t + o(t) = 1 + i m t - \sigma t + o(t)$ as $t \downarrow 0$}
\end{equation} 
for some $\beta = i m - \sigma$ $\in \CC$,
with $m \in \RR$ and $\sigma \geq 0$.  Let~$\nu$ be
the corresponding probability measure.  Then
$$
T_0^n \nu \implies \mbox{\rm Cauchy}(m, \sigma),
$$
where~$T_0$ is the homogeneous analogue of the {\tt Quicksort} transformation~$T$
mapping distributions as follows (in obvious notation):
\begin{equation}
\label{homo}
T_0: Z \mapsto U Z + (1 - U) Z^*.
\end{equation}
\end{theorem}

\begin{proof}
Let $Z_1, Z_2, \ldots; U_1, U_2, \ldots$ be independent random variables,
with every $Z_i \sim \nu$ and every $U_j \sim \mbox{unif}(0, 1)$.  Then,
using the definition of~$T_0$ repeatedly,
$$
W_n := \sum_{i = 1}^{2^n} V^{(n)}_i Z_i \sim T_0^n \nu, \qquad n \geq 0,
$$
where we define the random variables $V^{(n)}_i$ as follows.
Using~$U_1$ in the obvious fashion, split the unit interval into intervals of
lengths~$U_1$ and $1 - U_1$.  Now using~$U_2$ and~$U_3$, split the first interval
into subintervals of lengths $U_1 U_2$ and $U_1 (1 - U_2)$ and the second interval
into subintervals of lengths $(1 - U_1) U_3$ and $(1 - U_1) (1 - U_3)$.  Continue in this
way (using $U_1, \ldots, U_{2^n - 1}$) until the unit interval has been divided
overall into~$2^n$ subintervals.
Call their lengths, from left to right, $V^{(n)}_1, \ldots, V^{(n)}_{2^n}$.

Let $L_n := \max(V^{(n)}_1, \ldots, V^{(n)}_{2^n})$.  We show that~$L_n$
converges in probability to~$0$ as $n \to \infty$.  Luckily, the complicated
dependence structure of the variables~$V^{(n)}_i$ does not come into play;
the only observation we need is that that each $V^{(n)}_i$ marginally has the
same distribution as $U_1 \cdots U_n$.  Indeed, abbreviate $V^{(n)}_1 $ as~$V_n$;
briefly put, we derive a Chernoff's bound for~$\ln (1 / V_n)$ and then simply use
subadditivity.  To spell things out, let $x > 0$ be fixed and let $t \geq 0$.
Then 
$$
\PP(V_n \geq e^{-x}) \leq e^{t x} \EE V^t_n
  = e^{t x} \prod_{j = 1}^n \EE U^t_j
  = e^{t x} (1 + t)^{-n}
  = \exp[- (n \ln(1 + t) - x t)].
$$
Choosing the optimal $t = \frac{n}{x} - 1$ (valid for $n \geq x$), this yields
$$
\PP(V_n \geq e^{-x}) \leq \exp[- (n \ln (n/x) - n + x)]
  = \exp[- (n (\ln n - \ln (e x)) + x)]
$$
and thus
$$
\PP(L_n \geq e^{-x}) \leq 2^n \exp[- (n (\ln n - \ln (e x)) + x)]
  = \exp[- (n (\ln n - \ln (2 e x)) + x)] \to 0
$$
as $n \to \infty$.

Since~$L_n$ converges in probability to~$0$, we can therefore choose $\epsilon_n \to 0$ so
that $\PP(L_n > \epsilon_n) \to 0$.  To prove the theorem, it then suffices to prove
$$
\Wt_n := {\bf 1}(L_n \leq \epsilon_n) W_n \implies \mbox{Cauchy}(m, \sigma).
$$
For this, we note that the characteristic function~$\phi_n$ of~$\Wt_n$ is given
for $t \in \RR$ by 
\begin{equation}
\label{phint}
\phi_n(t) = \PP(L_n > \epsilon_n) + \EE \left[ {\bf 1}(L_n \leq \epsilon_n)
\prod_{i = 1}^{2^n} \psi(V^{(n)}_i t) \right].
\end{equation}
We will show that $\phi_n(t)$ converges to $e^{\beta t} = e^{i m t - \sigma t}$
for each fixed $t \geq 0$, and [since, further, $\phi_n(- t) \equiv \overline{\phi_n(t)}$]
this will complete the proof of the lemma.

Indeed, we need only consider the second term in~\eqref{phint}.
For that, the calculus estimates outlined in the proof of the lemma
preceding Theorem~7.1.2 in~\cite{Chung} demonstrate that,
when $L_n \leq \epsilon_n$,
$$
\prod_{i = 1}^{2^n} \psi(V^{(n)}_i t) = (1 + D_n) e^{\beta t}
$$
for complex random variables $D_n$ (depending on our fixed choice of $t \geq 0$)
satisfying $|D_n| \leq \delta_n$ for a deterministic sequence
$\delta_n\ [\ \equiv \delta(\epsilon_n t)] \to 0$ [with $\delta(s) \to 0$ as $s \to 0$]. 
[Leaving out the error estimates, the argument is
\begin{eqnarray*}
\log \left[ \prod_{i = 1}^{2^n} \psi(V^{(n)}_i t) \right]
  &\approx& \sum_{i = 1}^{2^n} \left( \psi(V^{(n)}_i t) - 1 \right) \\
  &\approx& \sum_{i = 1}^{2^n} \beta V^{(n)}_i t = \beta t. \qquad ]
\end{eqnarray*}
It now follows easily that $\phi_n(t) \to e^{\beta t}$, as desired. 
\end{proof}

Both the next lemma and its immediate corollary (\refL{L:convo}) will be used
in our proof of \refT{T:char}.

\begin{lemma}
\label{L:diff}
Let $\nu_i \in \Fc$, $i = 1, 2$.
Suppose that $(Z_1, Z_2)$ is a coupling of~$\nu_1$ and~$\nu_2$
such that the characteristic function of $Z_1 - Z_2$
satisfies~\eqref{linear}.  Then there exists a coupling
$(\Zt_1, \Zt_2)$ of~$\nu_1$ and~$\nu_2$ such that $\Zt_1 - \Zt_2 \sim$
{\rm Cauchy}$(m, \sigma)$.
\end{lemma}

\begin{proof}
Extend~$T$ to a transformation~$T_2$ on the class~$\Mc_2$
of probability measures on~$\RR^2$ by mapping
the distribution $\xi \in \Mc_2$ of $(X, Y)$ to the distribution~$T_2 \xi$ of 
$$
(U X + (1 - U) X^* + g(U),\ U Y + (1 - U) Y^* + g(U)),
$$
where~$U$, $(X, Y)$, and~$(X^*, Y^*)$ are independent, with $(X, Y) \sim \xi$,
$(X^*, Y^*) \sim \xi$, and $U \sim \mbox{unif}(0, 1)$, and where~$g$ is given
by~\eqref{gdef}.
(Note that we use the {\em same\/} uniform~$U$ for the~$Y$s as for the~$X$s!)
Of course, $T_2$ maps the marginal distributions
$\xi_1(\cdot) = \xi(\cdot \times \RR)$ of~$X$
and $\xi_2(\cdot) = \xi(\RR \times \cdot)$ of~$Y$ into~$T \xi_1$ and~$T \xi_2$, respectively;
more importantly for our purposes, it maps the distribution, call it~$\hat{\xi}$, of $X - Y$
into the distribution~$T_0 \hat{\xi}$, with~$T_0$ defined at~\eqref{homo}.

Now let~$\nu \in \Mc_2$ have marginals $\nu_i$, $i = 1, 2$.  Then ($T_2^n \nu)_{n \geq
1}$ has constant marginals $(\nu_1, \nu_2)$ as $n$ varies and so is a tight sequence.  We then
can find a weakly convergent subsequence, say,
$$
T_2^{n_k} \nu \implies \nu^{\infty} \in \Mc_2;
$$
of course, the limit $\nu^{\infty}$ again has marginals~$\nu_i$, $i = 1, 2$.
Moreover,
$$
T_0^{n_k} \hat{\nu} = \widehat{T_2^{n_k} \nu} \implies \widehat{\nu^{\infty}}.
$$

But, by supposition, the characteristic function of~$\hat{\nu}$ satisfies~\eqref{linear},
so \refT{T:attract} implies that $\widehat{\nu^{\infty}}$ is Cauchy$(m, \sigma)$.
Thus $\nu^{\infty} \in \Mc_2$ supplies the desired coupling.
\end{proof}

\begin{lemma}
\label{L:convo}
Let $\nu_i \in \Fc$, $i = 1, 2$.
Suppose that $(Z_1, Z_2)$ is a coupling of~$\nu_1$ and~$\nu_2$
such that $Z_1 - Z_2$
has zero mean and finite variance.  Then $\nu_1 = \nu_2$.
\nopf
\end{lemma}

\subsection{The proof}
\label{finish}

\addtolength{\textheight}{+0.3in} 

We now complete the proof of \refT{T:char}.

\begin{proof}
As discussed in Section~\ref{intro}, it is simple to check that $\Cc \subseteq \Fc$
(and that the elements of~$\Cc$ are all distinct).

Conversely, given $\nu \in \Fc$, let $Z_1 \sim \nu_1 := \nu$ and
$Z_2 \sim \nu_2 := \mu$ be independent random variables
(on some probability space); recall that~$\mu$
is the limiting {\tt Quicksort} measure, with zero
mean and finite variance.
Write $\psi_i$, $i = 1, 2$, for the characteristic
functions corresponding respectively to $\nu_i$, $i = 1, 2$.
By \refL{L:1} (or see \refC{C:expand}), $\psi_1$ satisfies~\eqref{linear}
[for some $(m, \sigma)$].
Of course, $\psi_2$ satisfies~\eqref{linear} with~$\beta$ taken to
be~$0$, so the characteristic function $t \mapsto \psi_1(t) \psi_2(-t)$ of $Z_1 - Z_2$
satisfies~\eqref{linear} for the same $(m, \sigma)$ as for~$\psi_1$.
Applying \refL{L:diff}, there exists a coupling $(\Zt_1, \Zt_2)$
of~$\nu_1$ and~$\nu_2$ such that $C := \Zt_1 - \Zt_2 \sim$\ Cauchy$(m, \sigma)$.
Without loss of generality (by building a suitable product space), we may
assume the existence of a random variable $Y \sim \mu$ on the same probability space
as~$\Zt_1$ and~$\Zt_2$ such that~$Y$ and~$C$ are independent. 

We know that the distribution~$\nu_1$ of $\Zt_1 = \Zt_2 + C$ is a fixed point of~$T$.
But so is the distribution~$\nu'_1 \in \Cc$ of $Z := Y + C$.
By \refL{L:convo} applied to $(\Zt_1, Z)$, $\nu = \nu_1 = \nu'_1 \in \Cc$, as desired.
\end{proof}

\end{document}